\documentclass[11pt]{amsart}
\usepackage{amssymb, latexsym, amsmath, amsfonts, tikz}
\usepackage{amsmath}
\usepackage{graphicx,color}
\usepackage{epstopdf}
\usepackage{epsfig}
\usepackage{verbatim}

\newtheorem{theorem}{Theorem}[section]
\newtheorem{corollary}[theorem]{Corollary}
\newtheorem{lemma}[theorem]{Lemma}
\newtheorem{proposition}[theorem]{Proposition}
\theoremstyle{definition}

\theoremstyle{remark}

\numberwithin{equation}{section}

\def\R{\mathbb{R}} 
\def\C{\mathbb{C}} 

\def\b{\mathbb{B}} 

\def\rt{\rightarrow}

\newcommand{\B}{\mathbb B}

\newcommand{\D}{\Delta}

\begin{document}

\sloppy

\title[Holomorphic motions and complex geometry]{Holomorphic motions and \\ complex geometry}

\author[Herv\'e Gaussier \and Harish Seshadri]{Herv\'e Gaussier and Harish Seshadri}

\address{\begin{tabular}{lll}
Herv\'e Gaussier & & Harish Seshadri\\
Univ. Grenoble Alpes, IF, F-38000 Grenoble, France & & Department of Mathematics\\
CNRS, IF, F-38000 Grenoble, France & &  Indian Institute of Science\\
& & Bangalore 560012\\
& & India
\end{tabular}
}
\email{herve.haussier@univ-grenoble-alpes.fr \ harish@math.iisc.ernet.in} 

\begin{abstract}
We show that the graph of a holomorphic motion of the unit disc cannot be biholomorphic to a strongly pseudoconvex domain in $\C^n$.
\end{abstract}


\maketitle 

\section{Introduction and Main Result}
Let $B$ be a connected complex $(n-1)$-manifold with a basepoint $z_0 \in B$. A {\it holomorphic motion} of  the unit disc $\D \subset \C$ parametrized by $B$ is a continuous map $f: B \times \D \rt \mathbb C \mathbb P^1 = \C \cup \{\infty\} $ satisfying the following conditions:

(1) $f(z_0,w)=w$ for all $w \in \D$,

(2) the map $f(z,.): \D \rt \mathbb C \mathbb P^1$ is injective for each $z \in B$,

(3) the map  $f(.,w): B \rt \mathbb C \mathbb P^1$ is holomorphic for each $w \in \D$. \vskip 0,2cm

Holomorphic motions were introduced by R. M\~ane, P. Sad and D. Sullivan \cite{mss} and have been intensively studied since then (see, for instance, \cite{slo, chi, dou, ast-mar}). In this note we study the complex-analytic structure of the graph $D$ of $f$:
$$D = \{(z,f(z,w)), \ z \in B, w \in \D \} \subset  B \times \C \mathbb P^1.$$
Our main result is the following

\begin{theorem}\label{str-pcv}
The graph $D$ of a holomorphic motion of the unit disc cannot be biholomorphic to a strongly pseudoconvex domain in $\C^n$.
\end{theorem}

Denoting the unit ball in $\C^n$ by $\B^n$, Theorem~\ref{str-pcv} will be a consequence of the following

\begin{theorem}\label{ball}
Let $S \subset \C$ be a bounded domain. If $A({\mathbb C}^n)$ denotes the set of complex affine $(n-1)$-dimensional subspaces of $\C^n$, then there does not exist a map $f: S \rt A({\mathbb C}^n)$ satisfying the following conditions: \vskip 0,2cm

(1) For $t \in S$, if $W_t = f(t) \cap \B^n$, then $\B^n= \cup_{t \in S} W_t$. \vskip 0,2cm

(2) Either $W_t \cap W_s = \phi$ or $W_s=W_t$ for $s,t \in S$.   \vskip 0,2cm

(3) There is a holomorphic map $\pi: \B^n \rt \B^{n-1}$  such that $\pi: W_t \rt \B^{n-1}$ is bijective for all $t \in S$.

\end{theorem}

To derive Theorem \ref{str-pcv} from Theorem \ref{ball}, we use a rescaling argument based on a  recent result of K. T. Kim and L. Zhang \cite{kz}.
\vspace{2mm}

To put these results in context, we recall the following result of K. Liu \cite{liu} and V. Koziarz-N. Mok \cite{km} :

 \begin{theorem}{\rm (\cite{liu}, \cite{km})} \label{nai}
Let $n > m \ge 1$ and let $\Gamma_1 \subset SU(n,1), \ \Gamma_2 \subset SU(m,1)$ be torsion-free cocompact lattices.
Then there does not exist a holomorphic submersion from ${\b ^n}/{\Gamma_1}$ to ${\b^m}/{\Gamma_2}$.
\end{theorem}

This was proved for $n=2, m=1$ by K. Liu \cite{liu} and  for all $n > m \ge 1$ by V. Koziarz and N. Mok \cite{km}. This result is natural from  various points of view. In particular, it is related to the following well-known question in the study  of negatively curved Riemannian manifolds:   \vskip 0,2cm

{\it Let $f: M^n \rt N^m$ be a smooth fibre bundle where $M$ and $N$ are smooth compact manifolds of dimensions $n > m \ge 2$. Can $M$ admit a Riemannian metric with negative sectional curvature ?}\vskip 0,2cm

If the bundle above is trivial then Preissman's theorem implies that the answer to the above question is in the negative. Also, it is essential that $m \ge 2$: a theorem of
W. Thurston states that certain $3$-manifolds fibering over a circle admit hyperbolic metrics.  \vskip 0,2cm

Theorem \ref{ball} implies the Liu-Koziarz-Mok result when $n=m+1$ by the Bers-Griffiths uniformization theorem as explained later in the paper.
The compactness of the manifolds is essential in the question above and the result of Liu-Koziarz-Mok. In other words, cocompact group actions on universal covers are needed. Our point of view is that given the Bers-Griffiths theorem, the cocompact actions are not necessary. The proof we present involves some  elementary 
facts about the Kobayashi metric and  Riemannian geometric techniques.  \vskip 0,2cm

Note that an equivalent formulation of Theorem \ref{ball} is that the graph of a holomorphic motion cannot admit a complete K\"ahler metric with constant negative holomorphic sectional curvature. Hence the following question is natural: \vspace{2mm}

{\it Can the graph of a holomorphic motion of the unit disc admit a complete K\"ahler metric with variable negative sectional curvature?} \vspace{2mm}

A related question, mentioned to the authors by Benoit Claudon and Pierre Py, is: \vspace{2mm}

{\it Can the graph of a holomorphic motion of the unit disc be Gromov hyperbolic with respect to the Kobayashi metric ?} \vspace{2mm}

 The method in this paper appears to hold some promise for tackling these questions. In this connection, it is important to point out that metrics with weaker negative curvature conditions can exist on such domains: a result of S. K. Yeung \cite{yeu} asserts that the universal cover of a  Kodaira fibration, which is necessarily the graph of a holomorphic motion by the Bers-Griffiths theorem, admits complete K\"ahler metrics with negative holomorphic bisectional curvature.

\section{Proof of Theorem~\ref{ball}}

\subsection{The Kobayashi metric on $D$} \ \ \ Let

$\bullet$ $B$ be a connected complex $(n-1)$-manifold with a basepoint $z_0 \in B$,

$\bullet$ $f: B \times \D \rt \mathbb C  \mathbb P^1$ a holomorphic motion,

$\bullet$ $D = \{(z,f(z,w)) \ : \ z \in  B, w \in \D \} \subset  B \times \D$,

$\bullet$ $F: B \times \D \rt D$ be defined by $F(z,w)=(z,f(z,w))$,

$\bullet$  $\pi: D \rt B$ denote the first projection,

$\bullet$ for $p \in D$, let $S_p= \pi^{-1}(\pi(p))$,

$\bullet$ for $w \in \D$, let $$F_w: B \rt D \ \ {\rm be} \ \  F_w(z)=f(z,w)  \ \ \ \ {\rm and} \ \ \ \ \Sigma_w= F_w(B).$$
\vspace{2mm} 

Note that $\pi^{-1}(z) =F(\{z\}\times \D)$ for every $z \in B$.

\begin{lemma}\label{lk}
For every $w \in \D$, the map $F:B \times \{w\} \rt D$ is a holomorphic embedding which is totally geodesic for the Kobayashi metrics on $B$ and $D$. \vspace{2mm}
\end{lemma}

{\bf Proof:} Since
$$ d_{\D}^K(z_1,z_2) \ge d_D^K(F(z_1,w),F(z_2,w))\ge d_{\D}^K(\pi \circ F(z_1,w), \pi \circ F(z_2,w))= d_{\D}^K(z_1,z_2)$$
by the distance decreasing property of the Kobayashi metric the result follows.

\subsection{Proof of Theorem \ref{ball}}
Let $S \subset \C$ be a bounded domain and $f: S \rt A({\mathbb C}^n)$ be a map, where $A({\mathbb C}^n)$ is the set of complex affine $(n-1)$-dimensional subspaces of $\C^n$. Suppose that we have \vspace{2mm}

(1) a partition $\B^n= \cup_{t \in S} W_t$ where $W_t = f(t) \cap \B^n$ and  \vspace{2mm}

(2) a holomorphic map $\pi: \B^n \rt \B^{n-1}$  such that $\pi \vert_{W_t}: W_t \rt \B^{n-1}$ is bijective. \vspace{2mm}

Before stating the next lemma, we recall that the Kobayashi metric on $\B^n$ coincides with the Bergman metric and is, in particular, a $C^2$ Riemannian metric.

\begin{lemma}
(1) For every $t \in S$ and $z \in \B^{n-1},$ $W_t \cap \pi^{-1}(z)$ consists of a single point and the intersection is orthogonal.
\vspace{2mm}

(2) The fibers of $\pi$ are equidistant, i.e., for any $z_1,z_2 \in \B^{n-1}$ and $p \in \pi^{-1}(z_1)$ we have
$$d_D^K(p, \ \pi^{-1}(z_2)) =d_D^K(\pi^{-1}(z_1), \ \pi^{-1}(z_2))$$ \vspace{2mm}

\end{lemma}\label{gf}
{\bf Proof:}
Let $t \in S$ and $z_1 \in \B^{n-1}$. It is clear that $W_t \cap \pi^{-1}(z_1)$ is a singleton since $\pi \vert_{W_t}: W_t \rt \B^{n-1}$ is bijective. Fix $p \in \pi^{-1}(z_1)$. Since   $\B^n= \cup_{t \in S} W_t$, $p \in W_t$ for some $t \in S$. Let $\gamma:[0,1] \rt D$ be a unit-speed geodesic with $\gamma(0) =p$ and $\gamma'(0) \in T_p W_t$. By Lemma \ref{lk} we can assume that $\gamma ([0,1]) \subset W_t$. Let $\gamma(1)= q \in \pi^{-1}(z) \cap W_t$.

   We claim that $\gamma$ is the shortest geodesic between $\pi^{-1}(z_1)$ and $\pi^{-1}(z)$. This is because
$$ d_D^K(p, q)) \ge d_{\B^{n-1}}^K(z_0,z)= l(\gamma) \ge d_D^K(p,q)$$
for any $w_1,w_2 \in \D$. The equality above comes from the assumption that $\pi \vert_{W_t}: W_t \rt \B^{n-1}$ is bijective and hence an isometry. Since $k$ is a Riemannian metric, the first variation for arc-length implies that $\gamma$ meets $\pi^{-1}(z_1)$ and
 $\pi^{-1}(z)$ orthogonally. Hence $\gamma'(0)$ is orthogonal to $T_{p} \pi^{-1}(z_1)$. \hfill $\square$
 \vspace{3mm}

In what follows we use the following notation: for any $p \in D$
$$S_p :=\pi^{-1}(\pi(p)).$$
\begin{corollary}\label{op}
Let $\gamma:[0,L] \rt D$ be a geodesic with $\gamma(0)=p, \ \gamma'(0) \in (T_p S_p) ^{\perp}$. If $P_s$ denotes the parallel transport of $T_p S_p$ along $\gamma$, then
$$P_s = T_{\gamma(s)} S_{\gamma(s)}.$$
\end{corollary}

{\bf Proof:} Let $\{e_1,e_2,...,e_{2n}\}$ be an orthonormal basis of $T_pD$ with $e_1, \ e_2  \in T_p S_p$ and let $E_i(s)$ be the parallel translate of $e_i$, $i=1,...,2n$, along $\gamma$.
By Lemma \ref{lk} for $i \ge 2$,  $\gamma$ lies in $W_t$ for some $t \in S$.  By Lemma (1) of \ref{gf}, $e_i \in T_pW_t$ for $3 \le i \le 2n$. Since $W_t$ is totally geodesic, $E_i(s)$ is tangent to $W_t$ for all $s \in [0,L]$ and $3 \le i \le 2n$. Hence
$E_1(s), E_2(s) \in (T_{\gamma(s)} W_t)^\perp = T_{\gamma(s)}S_{\gamma(s)}$. \hfill $\square$.

\subsection{Distance between complex submanifolds}
\vspace{2mm}

Fix $z_0,z_1 \in \B^{n-1}$. Let $p_0 \in \pi^{-1}(z_0),  p_1 \in \pi^{-1}(z_1)$ be points satisfying
$$d^K_D(p_0,p_1) =d_D^K(\pi^{-1}(z_0), \ \pi^{-1}(z_1)).$$

Let $\gamma:[0,L] \rt D$ be the unit speed geodesic with $\gamma(0)=p_0, \ \gamma(L)=p_1$  which realizes the distance between $\pi^{-1}(z_0)$ and $\pi^{-1}(z_1)$. Since $ \gamma'(0)$ is orthogonal to $T_p \pi^{-1}(z_0)$, Lemma \ref{lk} implies that $\gamma'(0) \in T_pW_t$ and $\gamma([0,\infty)) \subset W_t$ for some $t \in S$. \vspace{2mm}

By considering the normal exponential map to $\pi^{-1}(z_0)$, we can find a unit normal vector field $X$  to $\pi^{-1}(z_0)$ in a neighbourhood $U$ (in $\pi^{-1}(z_0)$) of $p_0$ such that the holds: for any $q \in U$, the geodesic $\gamma_q: [0,L] \rightarrow D$ with $\gamma_q(0)=q, \gamma'_q(0)=X_q$ satisfies 
$\gamma_q(L) \in \pi^{-1}(z_1)$. Note that $X_{p_0}= \gamma'(0)$. \vspace{2mm}

Let $u \in T_{p_0} \pi^{-1}(z_0)$ and $\sigma:[-a,a] \rt \pi^{-1}(z_0) \subset D$  a curve with $\sigma(0) = p_0$ and
$\sigma'(0)=u$.
Define a geodesic variation $H:[0,L] \times [-a,a] \rt D$ of $\gamma$ by
\begin{equation}
H(s,t) \ = \ Exp_{\sigma(t)} (sX_{\sigma(t)}).
\end{equation}

Let $Y(s,t)= \frac {\partial H}{\partial t}  (s,t)$ be the variation vector field and, for each $t \in [-a,a]$, let $\gamma_t$ be the geodesic given by
$\gamma_t(s)= H(s,t)$. Let $T(s,t)= \gamma_t'(s) = \frac {\partial H}{\partial s}  (s,t)$. \vspace{2mm}

\begin{lemma}\label{hg}
For any $(s,t) \in [0,L] \times [-a,a]$ we have \vspace{2mm}

(i) $Y(s,t) \in T_{\gamma_t(s)}S_{\gamma_t(s)}$.   \vspace{2mm}

(ii) $Y'(s,t) := \nabla_T Y (s,t) \in  T_{\gamma_t(s)} S_{\gamma_t(s)}$.
\end{lemma}

{\bf Proof:} (i) This follows if we can show that for each $s \in [0,L]$  the curve $t \mapsto H(s,t)$ lies in a fiber of $\pi$. To see this consider the curves, for $t_1, t_2
\in (-a,a)$, $s \mapsto (\pi \circ \gamma_{t_0})(s)$ and $s \mapsto (\pi \circ \gamma_{t_1})(s)$. These curves are unit-speed geodesics in $ \B^{n-1}$ connecting $\pi(z_0)$ and $\pi(z_1)$
by Lemma \ref{lk}. Since $ \B^{n-1}$ has negative curvature, uniqueness of geodesics forces $ \pi(H(s,t_0)) = \pi(H(s,t_1))$.

(ii) We show this for $t=0$ for notational simplicity. Let  $\{e_1, \ e_2 \}$ be an orthonormal basis of $T_{p_0} S_{p_0}$. Let $E_1(s), \ E_2(s)$ be the parallel vector fields along $\gamma$ with
$E_i(0)=e_i$. We then have, by (i) and Corollary \ref{op},
$$Y(s,0)= f_1(s)E_1(s) +f_2(s)E_2(s)$$
for some functions $f_1,f_2:[0,L] \rt \R$. Hence
$$Y'(s,0)= f_1'(s)E_1(s) +f_2'(s)E_2(s)  \ \in \   T_{\gamma(s)} S_{\gamma(s)}. \ \ \ \  \square$$

\vspace{3mm}
 We continue to  denote $p=F(z_0,w_0)$ in what follows. We recall the notation and constructions above:

$\bullet$ $X$ denotes a local unit  normal vector field on $S_p$ such that the geodesics starting in the direction $X$ pass through the same fibers subsequently,

$\bullet$ $\gamma:[0, \infty) \rightarrow D$ is a geodesic with $\gamma(0)=p, \gamma'(0) =X_p$,

$\bullet$ for $u \in T_p \ S_p$, $H_u(s,t)=\gamma^u_t(s)$ denotes a geodesic variation  of $\gamma$ such that

(a) $\frac {\partial H_u}{\partial s} (0,t) = X_{H(0,t)}$

(b) the variation vector field $Y_u(s,t)=\frac{\partial H_u}{\partial t}(s,t)$ satisfies
$Y_u(0,0)=u$,

$\bullet$ $T_u(s,t)=\frac {\partial H_u}{\partial s} (s,t)$ denotes the tangent vector $(\gamma_t^u)'(s)$ to the variation geodesic $\gamma_t^u(s)$.
Note that $T_u(s,0)= \gamma'(s)$ for all $s$ and $T_u(0,t) = X_{H_u(0,t)}$ for all $t$. \vspace{2mm}

\begin{lemma}\label{hh}
For any $u \in T_pS_p$ and $s \in [0,\infty)$, we have

(i) $ \langle \nabla_{Y_u}Y_u (s,0), \gamma'(s) \rangle = -\frac {1}{2} ( \vert Y_u \vert^2)'(s,0).$

(ii) $\nabla_uX  \ \in \ T_pS_p.$
\end{lemma}
{\bf Proof:}
\begin{align}\notag
\langle \nabla_{Y_u}Y_u, \gamma' \rangle \ &=  \ Y_u \langle Y_u, \gamma' \rangle - \langle Y_u, \nabla _{Y_u} \gamma' \rangle \\ \notag
&= \ - \langle Y_u, \nabla_{\gamma'} Y_u \rangle \\ \notag
&= \ -\frac {1}{2} ( \vert Y_u \vert^2)' \notag
\end{align}
where we have used $\langle T_u, Y_u \rangle =0$ and  $\nabla_{Y_u}T_u=\nabla_{T_u}Y_u$. This proves (i). For (ii), $[Y_u,T_u]=0$ and (ii) of Lemma \ref{hg} imply that  $\nabla_uX =\nabla_{\gamma'(0)}Y_u  \ \in \ T_pS_p.$ \hfill $\square$
\vspace{2mm}

By (ii) above, the {\it shape operator}  $L: T_p S_p \rt T_p S_p$ of $S_p$ along the normal vector field  $X$ is given by
$$L(v)= \nabla_v X.$$
As $L$ is a symmetric operator we can find an orthonormal basis $\{e_1,e_2\}$ of $T_p S_p$ consisting of eigenvectors of $L$. Since $S_p$ is a minimal submanifold (being a  complex subvariety),
the corresponding eigenvalues are given by $\alpha, -\alpha$ for some $\alpha \ge 0$.  As before we denote the parallel transports of $e_1,e_2$ along $\gamma$ by $E_1(s),E_2(s)$. \vspace{2mm}

Next we observe that if $Y_1 (s):=Y_{e_1}(s,0), \  Y_2(s):=Y_{e_2}(s,0)$ are Jacobi fields constructed as earlier with $Y_1(0)=e_1, \ Y_2(0)=e_2$ then
$$Y_1'(0) = \nabla_{\gamma'(0)} Y_1 = \nabla_{e_1}X = L(e_1) = \alpha e_1.$$
Similarly
$$Y_2'(0) = -\alpha e_2.$$

Hence if $K_1:[0,L] \rt \R$ denotes the function
$$ K_1(s) = R(E_1(s), \gamma'(s), \gamma'(s), E_1(s))$$
and $f_1: [0,L] \rt $ is the solution to
$$ y''+K_1y=0, \ \ \ y(0)=1, \ \ y'(0)=\alpha$$
then
$$Y_1 =f_1E_1.$$
Similarly $Y_2=f_2E_2$ where $f_2$ satisfies $y''+K_2y=0$, $y(0)=1$, $y'(0)=-\alpha$ with $K_2(s) =R(E_2(s), \gamma'(s), \gamma'(s), E_2(s))$.

\subsection{The case of $\C {\mathbb H}^n$}
In case $D$ is biholomorphic to the unit ball in $\C^n$, the Kobayashi metric on $D$ has constant holomorphic sectional curvature $-1$ and the curvature tensor
has the property that
$$\langle R(X,Y)Y,X \rangle = - \frac{1}{4}$$
whenever $\{X,Y\}$ is an orthonormal pair spanning a totally real $2$-plane, i.e., whenever $\langle X,Y \rangle =\langle X, JY \rangle =0$. Hence
$f_1$ and $f_2$ satisfy
$$y'' - \frac{y}{4}=0.$$
It follows that
$$f_1(s)= \cosh (\frac{s}{2}) + 2\alpha \sinh (\frac{s}{2}), \ \ \ \ f_2(s)=  \cosh (\frac{s}{2}) - 2 \alpha \sinh (\frac{s}{2}).$$
{\bf Case 1:}  $\alpha \neq \frac{1}{2}$.

In this case, (i) of Lemma \ref{hh} implies that
$$ \langle \nabla_{Y_i}Y_i, \gamma' \rangle (s) = - \frac {1}{2} (f^2_i)'(s) <0 $$
for $i=1,2$ and $s$ large enough. On the other hand, Lemma \ref{hg} and the fact that $S_{\gamma(s)}$ is a minimal submanifold implies that
$$0= \sum _{i=1}^2 \langle \nabla_{E_i}E_i, \gamma' \rangle (s)=\sum _{i=1}^2 f_i^{-2} \langle \nabla_{Y_i}Y_i, \gamma' \rangle (s).$$
This contradiction completes the proof. \vspace{2mm}

{\bf Case 2:} $\alpha = \frac{1}{2}$ for all $p \in D$, all $q \in S_p$ and all $T_0 \in (T_q S_p)^\perp$.
In this case, one can check that the second fundamental form of $S_p$ in every normal direction is parallel. O'Neill's formula \cite{one} (Page 465, 2. of Corollary 1) for the curvature of a Riemannian submersion then shows that the sectional curvature of the $2$-plane $\{u,v\}$ is zero where $u \in T_pS_p$ and $v \in (T_pS_p)^\perp$.

\subsection{Holomorphic motions and Bers-Griffiths uniformization}
The following fundamental theorem allows us to deduce (when $n=m+1$) Theorem~\ref{nai} from Theorem \ref{ball}:
\begin{theorem}(\cite{ber}, \cite{gri})
Let $M$ and $N$ be compact complex manifolds and $\phi:M \rt N$ a holomorphic submersion. Suppose that $dim(M)=dim(N)+1$ and the fibers of $\phi$ are compact Riemann surfaces of genus $\ge 2$.  Then the universal cover of $M$ is biholomorphic to the graph of a holomorphic motion over $\tilde N$.
\end{theorem}
For a detailed account of holomorphic motions and uniformization we refer the reader to \cite{chi}.

\section{Proof of Theorem~\ref{str-pcv}.}

We assume, to get a contradiction, that the graph $D$ of some holomorphic motion is biholomorphic to some bounded strictly pseudoconvex domain $\Omega$. We denote by $\Phi$ a biholomorphism from $D$ to $\Omega$.

For every $\nu \geq 1$, let $F_\nu : B \rt D$ be the totally geodesic holomorphic embedding defined by $F_{\nu}(z) := F(z,1-\frac{1}{\nu})$.

For every $\nu$, the map $\Phi \circ F_{\nu}$ is holomorphic from $B$ to $\Omega$. Let $z_0 \in B$. We may assume, taking a subsequence if necessary, that $\lim_{\nu \rt \infty} z^{\nu}:=\Phi \circ F_{\nu}(z_0) = p \in \partial \Omega$.

According to~\cite[Theorem 4.1]{kz} it holds $\lim_{z \rt p}\sigma_{\Omega}(z) = 1$, where $\sigma_{\Omega}$ is the squeezing function of $\Omega$ (see Definition in \cite{kz}).
This means that for every $\nu \geq 1$ there exists a biholomorphism $\varphi_{\nu}$ from $\Omega$ to some strongly pseudoconvex domain $\Omega_{\nu}$ and there exists a sequence $(r_{\nu})_{\nu}$ with $\lim_{\nu \rt \infty}r_{\nu} = 1$ such that for every $\nu \geq 1$:

 \begin{equation}\label{kz-eq}
 \varphi_{\nu}(z^{\nu}) = 0 \ {\rm and} \ B(0,r_{\nu}) \subset \Omega_{\nu} \subset \mathbb B^n.
 \end{equation}

Here $B(0,r_{\nu})$ denotes the ball in $\mathbb C^n$ centered at the origin with radius $r_{\nu}$.

For every $\nu \geq 1$, let $\Sigma_{\nu}:=\{F_{\nu}(z),\ z \in B\}$ and let $\tilde{\Sigma}_0^{\nu}:=(\varphi_{\nu} \circ \Phi)(\Sigma_{\nu})$.

\begin{lemma}\label{tot-lem}
For every $\nu \geq 1$, the set $\tilde{\Sigma}_0^{\nu}$ is a totally geodesic complex submanifold of $\Omega_{\nu}$.
\end{lemma}
\noindent{\bf Proof of Lemma~\ref{tot-lem}.}
This follows from Lemma~\ref{lk} and the fact that biholomorphisms are isometries for the Kobayashi metric. \hfill $\square$

Moreover, we get:

\begin{lemma}\label{conv-lem}
 The sequence $(\tilde{\Sigma}_0^{\nu})_{\nu}$ converges, for the local Hausdorff convergence of sets, to some totally geodesic complex submanifold of $\mathbb B^n$.
\end{lemma}

\noindent{\bf Proof of Lemma~\ref{conv-lem}.} Let, for every $\nu \geq 1$, $\Psi_{\nu}:= \varphi_{\nu} \circ \Phi \circ F_{\nu}$. Then $\Psi_{\nu}$ is a holomorphic isometric embedding of $(B,d^K_B)$ into $(\Omega_{\nu},d^K_{\Omega_{\nu}})$ satisfying $\Psi_{\nu}(z_0) = 0$. Since for every $\nu \geq 1$ we have the inclusion $\Omega_\nu \subset \mathbb B^n$, the sequence $(\Psi_{\nu})_{\nu}$ is normal and extracting a subsequence if necessary, we may assume that $(\Psi_{\nu})_{\nu}$ converges, uniformly on compact subsets of $B$, to some holomorphic map $\Psi_{\infty}:B \rt \mathbb B^n$ satisfying $\Psi_{\infty}(z_0) = 0$. Finally, let $ 0 \in L \subset \subset \mathbb B^n$. Since $\Psi_{\nu}$ is an isometry for the Kobayashi distances,  there exists $K \subset \subset B$ such that for every $\nu \geq 1$ we get: $L \cap \tilde{\Sigma}_0^{\nu} \subset \Psi_{\nu}(K)$. Now the uniform cnvergence of $(\Psi_{\nu})_{\nu}$ on $K$ implies that the sets $\tilde{\Sigma}_0^{\nu}$ converge to $\tilde{\Sigma}_0^{\infty}$ for
the Hausdorff convergence on $L$.

Let $\tilde{\Sigma}_0^{\infty}:=\Psi_{\infty}(B)$ and let $z,z' \in B$. There exist $q_{\nu}, q'_{\nu} \in \tilde{\Sigma}_0^{\nu}$, converging respectively to $\Psi_{\infty}(z)$ and $\Psi_{\infty}(z')$ and we have by Lemma~\ref{tot-lem}:
$$
\begin{array}{lllll}
d^K_{\tilde{\Sigma}_0^{\infty}}(\Psi_{\infty}(z),\Psi_{\infty}(z')) & = & \lim_{\nu \rt \infty}d^K_{\tilde{\Sigma}_0^{\nu}}(\Psi_{\nu}(z),\Psi_{\nu}(z')) & = &
\lim_{\nu \rt \infty}d^K_{\Omega_{\nu}}(\Psi_{\nu}(z),\Psi_{\nu}(z'))\\
& & & = & d^K_{\mathbb B^n}(\Psi_{\infty}(z),\Psi_{\infty}(z')).
\end{array}
$$
\qed

\vskip 0,2cm
In particular, since totally geodesic complex submanifolds of $\mathbb B^n$, of complex dimension $(n-1)$, are intersections of $\mathbb B^n$ with complex affine subspaces of complex dimension $(n-1)$, we may assume that $\tilde{\Sigma}_0^{\infty}=\mathbb B^{n-1} \times \{0\}$.

\vskip 0,2cm
Let $\zeta \in \Delta$ and let $q:=(0,\zeta) \in \mathbb B^n$. Then $q \in \Omega_{\nu}$ for sufficiently large $\nu$ and there exists $(b_{\nu},\zeta_{\nu}) \in B \times \Delta$ such that $q=\varphi_{\nu} \circ \Phi \circ F_{\zeta_{\nu}}(b_{\nu})$. We set $\tilde{\Sigma}^{\nu}_q:= \varphi_{\nu} \circ \Phi \circ F_{\zeta_{\nu}}(B)$. We prove, exactly as for $\tilde{\Sigma}^{\nu}_0$, that $\tilde{\Sigma}^{\nu}_q$ is a totally geodesic complex submanifold of $\Omega_{\nu}$.

\begin{lemma}\label{bded-lem}
 The sequence $(b_{\nu})_{\nu}$ is relatively compact in $B$.
\end{lemma}

\noindent{\bf Proof of Lemma~\ref{bded-lem}.} Since $D$ is complete hyperbolic by assumption,  it follows from Lemma~\ref{lk} that $B$ is complete hyperbolic. Assume to get a contradiction that $(b_{\nu})_{\nu}$ is not relatively compact in $B$. We recall that $\pi:D \rt B$ is holomorphic. Hence we get for every sufficiently large $\nu$:
$$
d^K_D\left(F\left(z_0,1-\frac{1}{\nu}\right), F(b_{\nu},\zeta_{\nu})\right) \geq \lim_{\nu \rt \infty}d^K_B(z_0,b_{\nu}).
$$
Consequently, extracting a subsequence if necessary, we may assume that:
$$
\lim_{\nu \rt \infty}d^K_D\left(F\left(z_0,1-\frac{1}{\nu}\right), F\left(b_{\nu},\zeta_{\nu}\right)\right)=\infty.
$$
Hence
\begin{equation}\label{infty-eq}
 \lim_{\nu \rt \infty}d^K_{\Omega_{\nu}}\left(\varphi_{\nu} \circ \Phi\left(F\left(z_0,1-\frac{1}{\nu}\right)\right), \varphi_{\nu} \circ \Phi(F(b_{\nu},\zeta_{\nu}))\right)=\infty.
\end{equation}

However, since $\varphi_{\nu} \circ \Phi\left(F\left(z_0,1-\frac{1}{\nu}\right)\right)=0$ and $\varphi_{\nu} \circ \Phi(F(b_{\nu},\zeta_{\nu}))=q$ for every $\nu$ we get:
$$
\lim_{\nu \rt \infty}d^K_{\Omega_{\nu}}\left(\varphi_{\nu} \circ \Phi\left(F\left(z_0,1-\frac{1}{\nu}\right)\right), \varphi_{\nu} \circ \Phi(F(b_{\nu},\zeta_{\nu}))\right)=d^K_{\mathbb B^n}(0,q) < \infty.
$$
This contradicts Condition~(\ref{infty-eq}). \qed

\vskip 0,2cm
It follows now from Lemma~\ref{bded-lem} that we may extract from $(b_{\nu})_{\nu})$ a subsequence, still denoted $(b_{\nu})_{\nu}$, that converges to some point $b_{\infty} \in B$.
Hence, extracting a subsequence if necessary, we may assume that the sequence $(\varphi_{\nu} \circ \Phi \circ F_{\zeta_{\nu}})_{\nu}$ converges uniformly on compact subsets of $B$ to a holomorphic map $\Psi^q_{\infty}:B \rt \mathbb B^n$ satisfying $\Psi^q_{\infty}(b_{\infty})= q$. This implies that $\tilde{\Sigma}^{\nu}_q=\varphi_{\nu} \circ \Phi \circ F_{\zeta_{\nu}}(B)$ converges to $\tilde{\Sigma}^{\infty}_q:=\Psi^q_{\infty}(B)$ and $\tilde{\Sigma}^{\infty}_q$ is a totally geodesic complex submanifold of $\mathbb B^n$.

We finally prove
\begin{proposition}\label{tot-prop}
 For every $q \in \mathbb B^n \cap (\{0'\} \times \D)$ there exists a totally geodesic complex submanifold $\tilde{\Sigma}^{\infty}_q$ of $\mathbb B^n$ passing through $q$.
\end{proposition}

For every $\nu$, let $\pi_{\nu} : D \rightarrow \Sigma_{\nu}$ be given by
$$
\forall (z,\zeta) \in B \times \Delta,\ \pi_{\nu}(F(z,\zeta)) = F\left(z,1-\frac{1}{\nu}\right)
$$
and let
$$
\begin{array}{ccccc}
 \tilde{\pi}_{\nu} & : & \Omega_{\nu} & \rightarrow & \Sigma_0^{\nu}\\
  & & z & \mapsto & \varphi_{\nu} \circ \Phi \circ \pi_{\nu} \circ \Phi^{-1} \circ \varphi_{\nu}^{-1}
\end{array}.
$$
Since $\varphi_{\nu} \circ \Phi(F(z_0,1-\frac{1}{\nu})) = 0$ according to (\ref{kz-eq}), we have $\tilde{\pi}_{\nu}(0) = 0$ for every $\nu$. Hence we may extract from $(\tilde{\pi}_{\nu})_{\nu}$ a subsequence, still denoted $(\tilde{\pi}_{\nu})_{\nu}$, that converges to a holomorphic map $\tilde{\pi}_{\infty} : \mathbb B^n \rightarrow \mathbb B^{n-1} \times \{0\}$.

Moreover we have:
\begin{proposition}\label{bihol-prop}
 For every $q \in \{0\} \times \Delta$, the restriction of $\tilde{\pi}_{\infty}$ to $\tilde{\Sigma}^{\infty}_q$ is a biholomorphism from $\tilde{\Sigma}^{\infty}_q$ to $\mathbb B^{n-1} \times \{0\}$.
\end{proposition}

\noindent{\bf Proof of Proposition~\ref{bihol-prop}.} By the very definition of $\tilde{\pi}_{\nu}$, the restriction of $\tilde{\pi}_{\nu}$ to $\tilde{\Sigma}_q^{\nu}$ is a biholomorphism from $\tilde{\Sigma}_q^{\nu}$ to $\mathbb B^{n-1} \times \{0\}$ for every $\nu$. Moreover, $\tilde{\Sigma}_q^{\nu}$ converges to $\tilde{\Sigma}_q^{\infty}$ for the Hausdorff distance. Finally, we have for every $\nu \geq 1$:
$$
\pi_{\nu} \circ \Phi^{-1} \circ \varphi_{\nu}^{-1}(q) = \left(b_{\nu},1-\frac{1}{\nu}\right) = F_{\nu}(b_{\nu}).
$$
Since $\lim_{\nu \rt \infty}b_{\nu} = b_{\infty} \in B$ and since the sequence $(\varphi_{\nu} \circ \Phi \circ F_{\nu})_{\nu}$ converges, uniformly on compact subsets of $B$ to $\Psi_{\infty}$ (see the proof of Lemma~\ref{conv-lem}) we obtain that:
$$
\tilde{\pi}_{\infty}(q) = \lim_{\nu \rt \infty}\varphi_{\nu} \circ \Phi \circ \pi_{\nu} \circ \Phi^{-1} \circ \varphi_{\nu}^{-1}(q) = \Psi_{\infty}(b_{\infty}) \in \mathbb B^{n-1} \times \{0\}.
$$
Hence if $g_{\nu}$ denotes the inverse of the restriction of $\tilde{\pi}_{\nu}$ to $\tilde{\Sigma}^{\infty}_q$ then $g_{\nu}$ is defined on $\mathbb B^{n-1} \times \{0\}$ and the sequence ($g_{\nu})_{\nu}$ converges, uniformly on compact subsets of $\mathbb B^{n-1} \times \{0\}$, to some holomorphic map $g_{\infty} : \mathbb B^{n-1} \times \{0\} \rt \tilde{\Sigma}^{\infty}_q$ such that $g_{\infty} \circ \tilde{\pi}_{\infty} = id_{\tilde{\Sigma}^{\infty}_q}$ and $\tilde{\pi}_{\infty} \circ g_{\infty}= id_{|\mathbb B^{n-1} \times \{0\}}$. \hfill $\square$

\vskip 0,2cm
Finally, let $q \neq q'$ be two points in $\mathbb B^n$. By construction, for every $\nu \geq 1$,  the intersection between the totally geodesics submanifolds $\tilde{\Sigma}_q^\nu$ and $\tilde{\Sigma}_{q'}^\nu$ is empty. Since $(\tilde{\Sigma}_q^\nu)_{\nu}$ converges to $\tilde{\Sigma}_q^{\infty}$ and $\tilde{\Sigma}_{q'}^\nu$ converges to $\tilde{\Sigma}_{q'}^{\infty}$, it follows from the positivity of intersection that:
$$
\tilde{\Sigma}_q^{\infty} \cap \tilde{\Sigma}_{q'}^{\infty} = \emptyset.
$$

Now Proposition~\ref{tot-prop} and Proposition~\ref{bihol-prop} give a contradiction, according to Theorem~\ref{ball}. \qed

\vskip 0,2cm
We end the note by studying some metric properties of $D$. We assume that $B$ is complete hyperbolic and that $B$ admits an exhaustion $(B_k)_{k \in \mathbb N}$: $B_{k} \subset \subset B_{k+1}$ for every $k$ and $B=\sup_{k \in \mathbb N} B_k$, such that $B_k$ is complete (Kobayashi) hyperbolic for every $k$.

\begin{proposition}\label{kob-prop}
The domain $D$ is complete (Kobayashi) hyperbolic.
\end{proposition}

\noindent{\bf Proof of Proposition~\ref{kob-prop}.} It is proved in \cite{che} that for every $k$, $D_k:=F(B_k \times \Delta)$ is complete hyperbolic.

Let $Z^0=(z_0,w_0), \ Z^{\nu}=(z_{\nu},w_{\nu}) \in B \times \Delta$ be such that $\lim_{\nu \rt \infty}d^K_B(z_0,z_{\nu}) = \infty$. Then:
\begin{equation}\label{inf-eq}
\forall \nu \geq 1,\ d^K_D(F(Z^0),F(Z^{\nu})) = d^K_B(z_0,z_{\nu}) \longrightarrow_{\nu \rt \infty} \infty.
\end{equation}

Hence, to prove that $D$ is complete hyperbolic, it is sufficient to prove that if $(z_{\nu})_{\nu} \subset \subset B_{k_0}$ for some $k_0 \in \mathbb N$ and $|w_{\nu}| \longrightarrow_{\nu \rt \infty}1$, then $d^K_D(F(Z^0),F(Z^{\nu})) \longrightarrow_{\nu \rt \infty} \infty$.
Assume, to get a contradiction, that there exists $c > 0$ such that $d^K_D(F(Z^0),F(Z^{\nu})) \leq c$ for every $k \in \mathbb N$ (extracting a subsequence if necessary). There exists $k_1 \geq k_0$ such that the set $\{y \in D /\ d^K_D(y,F(Z^0)) < c +1\}$ is contained in $D_{k_1}$ according to (\ref{inf-eq}).
Moreover, it follows from Lamma~5.1 in \cite{ki-ma} that:
$$
d^K_{D_{k_1}}(F(Z^0),F(Z^{\nu})) \leq \frac{1}{\tanh (1)}d^K_D(F(Z^0),F(Z^{\nu})).
$$
 This contradicts the fact that $D_{k_1}$ is complete hyperbolic. \qed

\end{document}